\documentclass[a4paper, 12pt]{article}
\usepackage{}

\usepackage{mathrsfs}
\usepackage{epsfig}
\usepackage{amsmath}
\usepackage{amssymb}
\usepackage{latexsym}
\usepackage{amsfonts}
\usepackage{amsthm}

\usepackage[numbers,sort&compress]{natbib}
\usepackage{graphicx}
\ifx\pdfoutput\undefined  \else
 \fi

\usepackage[centerlast]{caption2}

\usepackage{color}

\usepackage{txfonts}
\usepackage{amssymb}
\usepackage{mathrsfs}
\usepackage{amsfonts}

\newtheorem{theorem}{Theorem}[section]
\newtheorem{lemma}[theorem]{Lemma}

\usepackage{latexsym,amssymb}
\usepackage{amsmath}

\newcommand{\F}{{\mathbb F}}

\begin{document}

\title{A problem of Wang on Davenport constant for the multiplicative semigroup of the quotient ring of $\F_2[x]$}

\author{
Lizhen Zhang$^{a,b}$\thanks{Email:
lzhzhang0@aliyun.com}  \ \ \ \ \ \
Haoli Wang$^{c}$\thanks{Corresponding author's Email:
bjpeuwanghaoli@163.com}  \ \ \ \ \ \
Yongke Qu$^{d}$\thanks{Email:
yongke1239@163.com}
\\
\\
$^{a}${\small Shanghai Institute of Applied Mathematics and Mechanics} \\
{\small Shanghai University, Shanghai, 200072, P. R. China}\\
$^{b}${\small Department of Mathematics, Tianjin Polytechnic University, Tianjin, 300387, P. R. China}\\
$^{c}$ {\small College of Computer and Information Engineering}\\
{\small Tianjin Normal University, Tianjin, 300387, P. R. China}\\
$^{d}${\small Department of Mathematics, Luoyang Normal University,
Luoyang, 471022, P. R. China}
}

\date{}
\maketitle

\begin{abstract}  Let $\F_q[x]$ be the ring of polynomials over the finite field $\F_q$, and let $f$ be a polynomial of $\F_q[x]$.
Let $R=\frac{\F_q[x]}{(f)}$ be a quotient ring of $\F_q[x]$ with $0\neq R\neq \F_q[x]$. Let $\mathcal{S}_R$ be the multiplicative semigroup of the ring $R$, and let  ${\rm U}(\mathcal{S}_R)$ be the group of units of $\mathcal{S}_R$.
The Davenport constant ${\rm D}(\mathcal{S}_R)$ of the multiplicative semigroup $\mathcal{S}_R$ is the least positive integer $\ell$ such that for any $\ell$ polynomials $g_1,g_2,\ldots,g_{\ell}\in \F_q[x]$, there exists a subset $I\subsetneq [1,\ell]$  with
$$\prod\limits_{i\in I} g_i \equiv \prod\limits_{i=1}^{\ell} g_i\pmod f.$$
In this manuscript, we proved that for the case of $q=2$, $${\rm D}({\rm U}(\mathcal{S}_R))\leq {\rm D}(\mathcal{S}_R)\leq {\rm D}({\rm U}(\mathcal{S}_R))+\delta_f,$$
where
\begin{displaymath}
\delta_f=\left\{ \begin{array}{ll}
0 & \textrm{if $\gcd(x*(x+1_{\mathbb{F}_2}),\ f)=1_{\F_{2}}$}\\
1 & \textrm{if $\gcd(x*(x+1_{\mathbb{F}_2}),\ f)\in \{x, \ x+1_{\mathbb{F}_2}\}$}\\
2 & \textrm{if $gcd(x*(x+1_{\mathbb{F}_2}),f)=x*(x+1_{\mathbb{F}_2}) $}\\
\end{array} \right.
\end{displaymath}
which partially answered an open problem of Wang on Davenport constant for the multiplicative semigroup of $\frac{\F_q[x]}{(f)}$ (G.Q. Wang, \emph{Davenport constant for semigroups II,}
Journal of Number Theory, 155 (2015)  124--134).

\end{abstract}


\section {Introduction}

The additive properties of sequences in abelian groups have been widely studied in the field of Zero-sum Theory (see \cite{GaoGeroldingersurvey} for a survey), since H. Davenport \cite{Davenport} in 1966 and  K. Rogers \cite{rog1} in 1963 independently proposed one combinatorial invariant, denoted ${\rm D}(G)$, for any finite abelian group $G$, which is defined as the smallest $\ell\in \mathbb{N}$ such that
every sequence $T$ of terms from the group $G$ of length at least $\ell$ contains a nonempty subsequence $T'$ with sum of all terms from $T'$ being equal to the identity element of the group $G$. The Davenport constant is a central concept of zero-sum theory and has been investigated by many researchers in the scope of finite abelian groups.

In 2008, Gao and Wang \cite{wanggao} formulated the definition of Davenport constant for commutative semigroups, and made several related additive researches (see \cite{AdhikariGaoWang14, Wang II,  wang, wangirrecucibleseque}).

\noindent \textbf{Definition A.} \cite{wanggao} \ {\sl Let $\mathcal{S}$ be a commutative semigroup (not necessary finite). Let $T$ be a sequence of terms from the semigroup $\mathcal{S}$. We call $T$  reducible if $T$ contains a proper subsequence $T'$ ($T'\neq T$) such that the sum of all terms of $T'$ equals the sum of all terms of $T$. Define the Davenport constant of the semigroup $\mathcal{S}$, denoted ${\rm D}(\mathcal{S})$, to be the smallest $\ell\in \mathbb{N}\cup\{\infty\}$ such that every sequence $T$ of length at least $\ell$ of terms from $\mathcal{S}$ is reducible.}

Before then, starting from the research of Factorization Theory in Algebra, A. Geroldinger and F. Halter-Koch in 2006 have formulated another closely related definition, ${\rm d}(\mathcal{S})$, for any commutative semigroup  $\mathcal{S}$, which is called the small Davenport constant.

\noindent \textbf{Definition B.} (Definition 2.8.12 in \cite{GH}) \ {\sl For a commutative semigroup $\mathcal{S}$, let ${\rm d}(\mathcal{S})$ denote the smallest $\ell \in
\mathbb{N}_0\cup \{\infty\}$ with the following property:

For any $m\in \mathbb{N}$ and $a_1, \ldots,a_m\in \mathcal{S}$ there
exists a subset $I\subset [1,m]$ such that $|I|\leq \ell$ and
$$
\sum_{i=1}^m a_i=\sum_{i \in I}a_i.
$$}

The following connection between the (large) Davenport constant ${\rm D}(\mathcal{S})$ and the small Davenport constant ${\rm d}(\mathcal{S})$ was also obtained for any commutative semigroup $\mathcal{S}$.

\noindent \textbf{Proposition C.} (\cite{wangirrecucibleseque}) \ {\sl Let $\mathcal{S}$ be a commutative semigroup. Then ${\rm D}(\mathcal{S})$ is finite if and only if ${\rm d}(\mathcal{S})$ is finite. Moreover, in case that ${\rm D}(\mathcal{S})$ is finite, we have $${\rm D}(\mathcal{S})={\rm d}(\mathcal{S})+1.$$ }

Very recently, Wang in 2015 obtained the following result on Davenport constant for the multiplicative semigroup associated with polynomial rings $\F_q[x]$.

\noindent \textbf{Proposition D.} (\cite{Wang II}) \ {\sl Let $q>2$ be a prime power, and let $\F_q[x]$ be the ring of polynomials over the finite field $\F_q$.
Let $R$ be a quotient ring of $\F_q[x]$ with $0\neq R\neq \F_q[x]$. Then $${\rm D}(\mathcal{S}_R)={\rm D}({\rm U}(\mathcal{S}_R)),$$
where $\mathcal{S}_R$ denotes the multiplicative semigroup of the ring $R$, and ${\rm U}(\mathcal{S}_R)$ denotes the group of units in $\mathcal{S}_R$.}

However, for the case of $q=2$, Wang proposed it as an open problem.

\noindent \textbf{Problem E.} (see concluding remarks in \cite{Wang II}) \ {\sl Let $R$ be a quotient ring of $\mathbb{F}_2[x]$ with $0\neq R\neq \mathbb{F}_2[x]$. Determine ${\rm D}(\mathcal{S}_R)-{\rm D}({\rm U}(\mathcal{S}_R))$.}

In this manuscript, we considered this open problem. By using the method employed by Wang, we obtained the following result, which is a partial solution of Problem E.

\bigskip

\begin{theorem}\label{Theorem main}
\  Let $\F_2[x]$ be the ring of polynomials over the finite field $\F_2$, and let $R=\frac{\F_2[x]}{(f)}$ be a quotient ring of $\F_2[x]$ where $f\in\F_2[x]$ and $0\neq R\neq \F_2[x]$.
Then  $${\rm D}({\rm U}(\mathcal{S}_R))\leq {\rm D}(\mathcal{S}_R)\leq {\rm D}({\rm U}(\mathcal{S}_R))+\delta_f,$$
where
\begin{displaymath}
\delta_f=\left\{ \begin{array}{ll}
0 & \textrm{if $\gcd(x*(x+1_{\mathbb{F}_2}),\ f)=1_{\F_{2}}$;}\\
1 & \textrm{if $\gcd(x*(x+1_{\mathbb{F}_2}),\ f)\in \{x, \ x+1_{\mathbb{F}_2}\}$;}\\
2 & \textrm{if $gcd(x*(x+1_{\mathbb{F}_2}),f)=x*(x+1_{\mathbb{F}_2})$.}\\
\end{array} \right.
\end{displaymath}
\end{theorem}

\section{The proof of Theorem \ref{Theorem main}}

The notations and terminologies used here are consistent to ones used in \cite{{AdhikariGaoWang14}, Wang II, wangirrecucibleseque, wang}. For the reader's convenience, we need to give some necessary ones.

Let $\mathcal{S}$ be a finite commutative semigroup.
The operation on $\mathcal{S}$ is denoted by $+$.
The identity element of $\mathcal{S}$, denoted $0_{\mathcal{S}}$ (if exists), is the unique element $e$ of
$\mathcal{S}$ such that $e+a=a$ for every $a\in \mathcal{S}$. If $\mathcal{S}$ has an identity element $0_{\mathcal{S}}$, let
$${\rm U}(\mathcal{S})=\{a\in \mathcal{S}: a+a'=0_{\mathcal{S}} \mbox{ for some }a'\in \mathcal{S}\}$$ be the group of units
of $\mathcal{S}$. For any element $c\in\mathcal{S}$, let $${\rm St}(c)=\{a\in {\rm U}(\mathcal{S}): a+c=c\}$$ denote the stabilizer of $c$ in the group ${\rm U}(\mathcal{S})$.
The Green's preorder of the semigroup $\mathcal{S}$, denoted $\leqq_{\mathcal{H}}$, is defined by
$$a \leqq_{\mathcal{H}} b\Leftrightarrow a=b \ \ \mbox{or}\ \ a=b+c$$ for some $c\in \mathcal{S}$. The Green's congruence of $\mathcal{S}$, denoted
$\mathcal{H}$, is defined by:
$$a \ \mathcal{H} \ b \Leftrightarrow a \ \leqq_{\mathcal{H}} \ b \mbox{ and } b \ \leqq_{\mathcal{H}} \ a.$$
We write $a<_{\mathcal{H}} b$ to mean that $a \leqq_{\mathcal{H}} b$ but $a \ \mathcal{H} \ b$ does not hold.

The sequence $T$ of terms from the semigroups $\mathcal{S}$ is denoted by $$T=a_1a_2\cdot\ldots\cdot a_{\ell}=\coprod\limits_{a\in \mathcal{S}}a^{[{\rm v}_a(T)]},$$ where $[{\rm v}_a(T)]$ means that the element $a$ occurs ${\rm v}_a(T)$ times in the sequence $T$.
By $\cdot$ we denote the operation to join sequences.
By $|T|$ we denote the length of the sequence, i.e., $$|T|=\sum\limits_{a\in \mathcal{S}}{\rm v}_a(T)=\ell.$$
Let $T_1,T_2$ be two sequences of terms from the semigroups $\mathcal{S}$. We call $T_2$
a subsequence of $T_1$ if $${\rm v}_a(T_2)\leq {\rm v}_a(T_1)$$ for every element $a\in \mathcal{S}$, denoted by $$T_2\mid T_1.$$ In particular, if $T_2\neq T_1$, we call $T_2$ a {\sl proper} subsequence of $T_1$, and write $$T_3=T_1  T_2^{[-1]}$$ to mean the unique subsequence of $T_1$ with $T_2\cdot T_3=T_1$.  Let $$\sigma(T)=a_1+a_2+\cdots+a_{\ell}$$ be the sum of all terms of the sequence $T$.
By $\varepsilon$ we denote the
empty sequence.
If $S$ has an identity element $0_{\mathcal{S}}$,  we allow $T=\varepsilon$ and adopt the convention
that $\sigma(\varepsilon)=0_\mathcal{S}$.
We say that $T$ is {\it
reducible} if $\sigma(T')=\sigma(T)$ for some proper subsequence $T'$ of $T$
(note that, $T'$ is probably the empty sequence $\varepsilon$ if $\mathcal{S}$
has the identity element $0_{\mathcal{S}}$ and $\sigma(T)=0_{\mathcal{S}}$). Otherwise, we call $T$
{\it irreducible}.

\bigskip

Throughout this paper, we shall always denote $$R=\mathbb{F}_2[x]\diagup (f)$$ to be the quotient ring  of $\mathbb{F}_2[x]$ modulo some nonconstant  polynomial $f\in \mathbb{F}_2[x]$, where
\begin{equation}\label{equation factorization of f(x)}
f=f_1^{n_{1}}*f_2^{n_{2}}*\cdots * f_r^{n_r},
\end{equation}
 such that $f_1,f_2, \ldots,f_r$ are pairwise non-associate irreducible polynomials of $\mathbb{F}_2[x]$ with $$f_1=x, \ f_2=x+1_{\mathbb{F}_2},$$ $$n_{1}\geq 0,n_2\geq 0, n_3,n_4,\ldots,n_r\geq 1.$$ Let $\mathcal{S}_R$ be the multiplicative semigroup of the ring $R$. For any element $a\in \mathcal{S}_R$, let $\theta_a\in \mathbb{F}_2[x]$ be the unique polynomial corresponding to the element $a$ with the least degree, i.e.,
$$\overline{\theta_a}=\theta_a+(f)$$ is the corresponding form of $a$ in the quotient ring $R$ with $${\rm deg}(\theta_a)\leq {\rm deg}(f)-1.$$
By $\gcd(\theta_a,f)$ we denote the greatest common divisor  of the two polynomials $\theta_a$ and $f$ in $\mathbb{F}_2[x]$ ({\sl the unique polynomial with the greatest degree which divides both $\theta_a$ and $f$}), in particular, by \eqref{equation factorization of f(x)}, we have
\begin{equation}\label{equation factorization of gcdfa}
\gcd(\theta_a,f)=f_1^{\alpha_{1}}\ast f_2^{\alpha_{2}} *\cdots * f_r^{\alpha_r}
\end{equation}
where
$\alpha_i\in [0,n_i]$ for each $i\in [1,r]$.

For any polynomial $g$ and any irreducible polynomial $h$ of $\mathbb{F}_2[x]$, let ${\rm pot}_{h}(g)$ be the largest integer $k$ such that $h^k \mid g.$ Then in \eqref{equation factorization of gcdfa}, $\alpha_i={\rm pot}_{f_i}(\gcd(\theta_a,f))$ for each $i\in [1,r]$.
It is easy to observe that for any two element $a,b\in \mathcal{S}_R$, $$\gcd(\theta_a,f)=\gcd(\theta_b,f)$$ if and only if $${\rm pot}_{f_i}(\gcd(\theta_a,f))={\rm pot}_{f_i}(\gcd(\theta_b,f))\mbox{ for each }i\in [1,r].$$

To prove Theorem \ref{Theorem main}, we still need some lemmas.

\medskip

\begin{lemma}(\cite{GH}, Lemma 6.1.3) \label{Lemma recusive Davenport constant} \ Let $G$ be a finite abelian group, and let $H$ be a subgroup of $G$. Then, ${\rm D}(G)\geq {\rm D}(G/H)+{\rm D}(H)-1$.
\end{lemma}

\begin{lemma}\label{proposition D(U(G))leq D(G)} (see \cite{wanggao}, Proposition 1.2) \
Let $\mathcal{S}$ be a finite commutative semigroup with an identity. Then ${\rm D}({\rm U}(\mathcal{S}))\leq
{\rm D}(\mathcal{S})$.
\end{lemma}

\begin{lemma}\label{stabilizer1} \ Let $a$ and $b$ be two elements of $\mathcal{S}_R$ with $a \leqq_{\mathcal{H}} b$. Let $\alpha_i={\rm pot}_{f_i}(\gcd(\theta_a,f))$ and $\beta_i={\rm pot}_{f_i}(\gcd(\theta_b,f))$ for each $i\in [1,r]$. Then,

(i). ${\rm St}(b)\subseteq{\rm St}(a)$ and $\beta_i\leq \alpha_i$ for each $i\in [1,r]$, in particular, if $a\ \mathcal{H}\ b$ then ${\rm St}(b)={\rm St}(a)$ and $\beta_i=\alpha_i$ for each $i\in [1,r]$;

(ii).  if $\beta_i=\alpha_i$ for each $i\in [1,r]$,  then $a\ \mathcal{H}\ b$;

(iii). If $a<_{\mathcal{H}} b$ and $(\alpha_1-\beta_1)(2n_1-1-\alpha_1-\beta_1)+(\alpha_2-\beta_2)(2n_2-1-\alpha_2-\beta_2)+\sum^{r}_{i=3 }(\alpha_i-\beta_i)>0,$
then ${\rm St}(b)\subsetneq {\rm St}(a).$
\end{lemma}

\noindent {\sl Proof of Lemma \ref{stabilizer1}.} \
Note first that $a \leqq_{\mathcal{H}} b$ implies that $$\alpha_i\geq \beta_i \mbox{ for each }i\in [1,r].$$

(i). Since $\mathcal{S}_R$ has the identity element $0_{\mathcal{S}_R}$, it follows from $a \leqq_{\mathcal{H}} b$ that $$a=b+c \ \ \mbox{for some }c\in \mathcal{S}_R.$$ It follows that
$$\gcd(\theta_b,f)\mid \gcd(\theta_b*\theta_c,f)=\gcd(\theta_a,f),$$ equivalently, $\beta_i\leq \alpha_i$ for each $i\in [1,r]$.

Take an arbitrary element $d\in {\rm St}(b)$. Then  $d+a=d+(b+c)=(d+b)+c=b+c=a$,
and so $d\in {\rm St}(a)$. It follows that $${\rm St}(b)\subseteq {\rm St}(a).$$

If $a\ \mathcal{H}\ b$, i.e., $a \leqq_{\mathcal{H}} b$ and $b \leqq_{\mathcal{H}} a$, then
${\rm St}(b)={\rm St}(a)$ and $\beta_i=\alpha_i$ for each $i\in [1,r]$ follows readily.
This proves Conclusion {\rm (i)}.

(ii). Assume $\beta_i=\alpha_i$ for each $i\in [1,r]$, that is,  $$\gcd(\theta_b,f)=\gcd(\theta_a,f).$$
It follows that there exist polynomials $h, h'\in \mathbb{F}_q[x]$ such that $$\theta_a* h\equiv \theta_b\pmod {f}$$
and
$$\theta_b* h'\equiv \theta_a\pmod {f}.$$
It follows that $b \leqq_{\mathcal{H}}  a$ and $a \leqq_{\mathcal{H}}  b$, i.e., $$a \ \mathcal{H} \ b,$$ and Conclusion {\rm (ii)} is proved.

(iii). Now assume $$a<_{\mathcal{H}} b$$ and
\begin{equation}\label{equation three sum}
\sum^{r}_{i=3 }(\alpha_i-\beta_i)+(\alpha_1-\beta_1)(2n_1-1-\alpha_1-\beta_1)+(\alpha_2-\beta_2)(2n_2-1-\alpha_2-\beta_2) >0.\end{equation}  It is sufficient to find some element $d\in {\rm U}(\mathcal{S}_R)$ such that $d\in {\rm St}(a)\setminus {\rm St}(b).$
We shall distinguish two cases.

\textbf{Case 1.} $\sum^{r}_{i=3 }(\alpha_i-\beta_i)>0$.

Then there exists some $i \in [3, r]$ such that $\alpha_{i}>\beta_{i}$, say
\begin{equation}\label{equation a3>b3}
\alpha_{3}>\beta_{3}.
\end{equation}
Take an polynomial
\begin{equation}\label{equation h=in case 1}
h=\frac{f}{f_3^{\alpha_3}}.
\end{equation}

We show that
\begin{equation}\label{equation h_2(x)+1,f(x)=1}
\gcd(h+1_{\mathbb{F}_2},f)=1_{\mathbb{F}_2}
\end{equation}
or
\begin{equation}\label{equation 2h_2(x)+1,f(x)=1}
 \gcd(x*h+1_{\mathbb{F}_2},f)=1_{\mathbb{F}_2}.
 \end{equation}

Suppose to the contrary that $\gcd(h+1_{\mathbb{F}_2},f)\neq 1_{\mathbb{F}_2}$ and $\gcd(x*h+1_{\mathbb{F}_2},f)\neq 1_{\mathbb{F}_2}$.
By \eqref{equation factorization of f(x)} and \eqref{equation h=in case 1}, we have that $f_i\not\mid \gcd(h+1_{\mathbb{F}_2},f)$ and $f_i\not\mid \gcd(x*h+1_{\mathbb{F}_2},f)$ for each $i\in[1,r]\setminus \{3\}$. This implies that $f_3\mid (h+1_{\mathbb{F}_2})$ and $f_3 \mid (x*h+1_{\mathbb{F}_2})$, and thus, $f_3\mid x*(h+1_{\mathbb{F}_2})-(x*h+1_{\mathbb{F}_2})=x+1_{\mathbb{F}_2}$, which is absurd. This proves that \eqref{equation h_2(x)+1,f(x)=1} or \eqref{equation 2h_2(x)+1,f(x)=1} holds.

Take an element $d\in \mathcal{S}_R$ with $$\theta_{d} \equiv h+1_{\mathbb{F}_{2}}\pmod {f}$$ or $$\theta_{d}\equiv x*h+1_{\mathbb{F}_2} \pmod {f}$$ according to \eqref{equation h_2(x)+1,f(x)=1} or \eqref{equation 2h_2(x)+1,f(x)=1} holds respectively. It follows that $$d\in {\rm U}(\mathcal{S}_R),$$ and follows from \eqref{equation a3>b3} and \eqref{equation h=in case 1} that
$$\theta_a*\theta_{d}\equiv \theta_a\pmod {f}$$
and $$\theta_b*\theta_{d}\not\equiv \theta_b\pmod {f}.$$ That is, $d\in {\rm St}(a)\setminus {\rm St}(b)$, which implies  $${\rm St}(b)\subsetneq {\rm St}(a).$$

\textbf{Case 2.} \ $(\alpha_1-\beta_1)(2n_1-1-\alpha_1-\beta_1)>0$ or $(\alpha_2-\beta_2)(2n_2-1-\alpha_2-\beta_2)>0$.

Say $$(\alpha_1-\beta_1)(2n_1-1-\alpha_1-\beta_1)>0.$$

It follows that
\begin{equation}\label{equation case two equ1}
\alpha_1>\beta_1
\end{equation}
and
\begin{equation}\label{equation case two equ2}
n_1>\beta_1+1.
\end{equation}
Take an polynomial \begin{equation}\label{equation case two equ3}h=\frac{f}{f_1^{\beta_1+1}}.
\end{equation}
Combined with \eqref{equation case two equ2} and \eqref{equation case two equ3}, we conclude that  $$\gcd(h+1_{\mathbb{F}_2},f)=1_{\mathbb{F}_2}.$$ Take an element $d\in \mathcal{S}_R$ with $$\theta_{d} \equiv h+1_{\mathbb{F}_{2}}\pmod {f}.$$ It follows that $$d\in {\rm U}(\mathcal{S}_R),$$ and follows from \eqref{equation case two equ1} and \eqref{equation case two equ3} that
$$\theta_a*\theta_{d}\equiv \theta_a\pmod {f}$$
and $$\theta_b*\theta_{d}\not\equiv \theta_b\pmod {f}.$$ That is, $d\in {\rm St}(a)\setminus {\rm St}(b)$ which implies  $${\rm St}(b)\subsetneq {\rm St}(a).$$
This proves Lemma \ref{stabilizer1}. \qed

\bigskip

Now we are in a position to prove Theorem \ref{Theorem main}.

\noindent {\sl Proof of Theorem \ref{Theorem main}.} \
By Lemma \ref{proposition D(U(G))leq D(G)}, it suffices to show that ${\rm D}(\mathcal{S}_R)\leq {\rm D}({\rm U}(\mathcal{S}_R))+\delta_f$.
Let $T=a_1a_2\cdot\ldots\cdot a_{\ell}$ be an arbitrary sequence of terms from $\mathcal{S}_R$ of length 
\begin{equation}\label{equation length of ell}
\ell={\rm D}({\rm U}(\mathcal{S}_R))+\delta_f.
 \end{equation}
 We shall prove that $T$ contains a {\sl proper} subsequence $T'$ with $\sigma(T')=\sigma(T)$.

Take a shortest subsequence $V$ of $T$ such that
\begin{equation}\label{equation sigma(V)Hsigma(T)}
\sigma(V) \ \mathcal{H} \ \sigma(T).
\end{equation}
We may assume without loss of generality that $$V=a_1\cdot a_2\cdot\ldots\cdot a_t\ \ \ \mbox{where} \ \ t\in [0,\ell].$$ By the minimality of $|V|$, we derive that $$0_{\mathcal{S}_R}>_{\mathcal{H}}a_1>_{\mathcal{H}}(a_1+a_2)>_{\mathcal{H}}
\cdots>_{\mathcal{H}}\sum_{i=1}^t a_i.$$
Denote $$K_0=\{0_{\mathcal{S}_R}\}$$
and
$$K_i={\rm St}(\sum\limits_{j=1}^i a_j) \ \ \ \mbox{for each}\ \ \ i\in [1,t].$$
Note that $K_i$ is a subgroup of ${\rm U}(\mathcal{S}_R)$ for each $i\in [1,t]$. Moreover, since ${\rm St}(0_{\mathcal{S}_R})=K_0$,
it follows from Conclusion (i) of Lemma \ref{stabilizer1} that
$$K_0\subseteq K_1\subseteq K_2\subseteq \cdots\subseteq K_t,$$
and moreover, by applying Lemma \ref{stabilizer1}, we conclude that
there exists a subset $M$ of $[0,t-1]$ with
\begin{equation}\label{equation t minus M geq}
|M|\geq t-\delta_f
 \end{equation}
 such that $$K_i\subsetneq K_{i+1}\mbox{ for each } i\in M.$$
For $i\in M$, since $\frac{{\rm U}(\mathcal{S}_R)}{K_{i+1}}\cong \frac{{\rm U}(\mathcal{S}_R)\diagup K_{i}}{K_{i+1}\diagup K_{i}}$ and ${\rm D}(K_{i+1}\diagup K_{i})\geq 2$, it follows from Lemma \ref{Lemma recusive Davenport constant} that
\begin{equation}\label{equation case for M}
\begin{array}{llll}
{\rm D}({\rm U}(\mathcal{S}_R)\diagup K_{i+1})&=& {\rm D}(\frac{{\rm U}(\mathcal{S}_R)\diagup K_{i}}{K_{i+1}\diagup K_{i}}) \\
&\leq & {\rm D}({\rm U}(\mathcal{S}_R)\diagup K_{i})-({\rm D}(K_{i+1}\diagup K_{i})-1)\\
&\leq & {\rm D}({\rm U}(\mathcal{S}_R)\diagup K_{i})-1.\\
\end{array}
\end{equation}
Combined with \eqref{equation length of ell}, \eqref{equation t minus M geq} and \eqref{equation case for M}, we conclude that
 \begin{align}\label{equation length and D()}
\begin{array}{llll}
1\leq {\rm D}({\rm U}(\mathcal{S}_R)\diagup K_t)
&\leq & {\rm D}({\rm U}(\mathcal{S}_R)\diagup K_0)-|M|\\
&\leq & {\rm D}({\rm U}(\mathcal{S}_R))-(t-\delta_f)\\
&=& (\ell-\delta_f)-(t-\delta_f)\\
&=& \ell-t\\
&=& |TV^{[-1]}|.\\
\end{array}
\end{align}
By Conclusion {\rm (i)} of Lemma \ref{stabilizer1} and \eqref{equation sigma(V)Hsigma(T)}, we have
\begin{equation}\label{equation two common divisors equal}
{\rm pot}_{f_i}(\gcd(\theta_{\sigma(V)},f))={\rm pot}_{f_i}(\gcd(\theta_{\sigma(T)},f))
\end{equation}
for each $i\in [1,r]$.
Let $$\mathcal{J}=\{j\in [1,r]: f_j^{n_j}\mid \theta_{\sigma(T)}\}.$$ By \eqref{equation two common divisors equal}, we have that
\begin{equation}\label{equation fi(x)notmid a}
f_i\not\mid\theta_{a} \ \ \ \mbox{for each term} \ a \ \mbox{of} \ TV^{[-1]} \ \mbox{and each} \ i\in [1,r]\setminus \mathcal{J},
\end{equation}
and that
\begin{equation}\label{equation fjnj(x)mid sigma(V)}
f_j^{n_j}\mid \theta_{\sigma(V)}  \ \ \ \mbox{ for each } j\in \mathcal{J}.
\end{equation}

For each term $a$ of $TV^{[-1]}$, let $\tilde{a}$ be the element of $\mathcal{S}_R$ such that
\begin{equation}\label{equation tilde a 1}
\theta_{\tilde{a}}\equiv \theta_{a}\pmod {f_i^{n_i}} \ \ \ \mbox{ for each } i\in [1,r]\setminus \mathcal{J}
\end{equation}
and
\begin{equation}\label{equation tilde a 2}
\theta_{\tilde{a}}\equiv 1_{\mathbb{F}_2}\pmod {f_j^{n_j}} \ \ \ \mbox{ for each } j \in\mathcal{J}.
\end{equation}
By \eqref{equation fi(x)notmid a}, \eqref{equation tilde a 1} and \eqref{equation tilde a 2}, we conclude that $\gcd(\theta_{\tilde{a}},f)=1_{\mathbb{F}_2}$, i.e.,
\begin{equation}\label{equation tilde a in U}
\tilde{a}\in {\rm U}(\mathcal{S}_R)\ \ \ \mbox{for each term} \ a \ \mbox{of}\  TV^{[-1]}.
\end{equation}
By \eqref{equation fjnj(x)mid sigma(V)} and \eqref{equation tilde a 1}, we conclude that \begin{equation}\label{equation sigma(V)+tilde a}
\sigma(V)+\tilde{a}=\sigma(V)+a \ \ \ \mbox{for each term}\ \  a \ \mbox{of} \ TV^{[-1]}.
\end{equation}
By \eqref{equation length and D()} and \eqref{equation tilde a in U}, we have that $\coprod\limits_{a\mid TV^{[-1]}}\tilde{a}$ is a nonempty sequence of elements in ${\rm U}(\mathcal{S}_R)$ of length $|\coprod\limits_{a\mid TV^{[-1]}}\tilde{a}|=|TV^{[-1]}|\geq {\rm D}({\rm U}(\mathcal{S}_R)\diagup K_t)$. It follows that there exists a {\bf nonempty} subsequence $$W\mid TV^{[-1]}$$  such that
$$\sigma(\coprod\limits_{a\mid W}\tilde{a})\in K_t$$ which implies
\begin{equation}\label{equation sigma(tilde a)=0}
\sigma(V)+\sigma(\coprod\limits_{a\mid W}\tilde{a})=\sigma(V).
\end{equation}
By \eqref{equation sigma(V)+tilde a} and \eqref{equation sigma(tilde a)=0}, we conclude that
$$\begin{array}{llll}
\sigma(T)&=& \sigma(TW^{[-1]}V^{[-1]})+(\sigma(V)+\sigma(W)) \\
&=& \sigma(TW^{[-1]}V^{[-1]})+(\sigma(V)+\sigma(\coprod\limits_{a\mid W}\tilde{a}))\\
&=& \sigma(TW^{[-1]}V^{[-1]})+\sigma(V)\\
&=& \sigma(TW^{[-1]}),\\
\end{array}$$
and $T'=TW^{[-1]}$ is the desired proper subsequence of $T$.
This completes the proof of the theorem. \qed

\bigskip

\noindent {\bf Acknowledgements}

\noindent 
This work is supported by NSFC (11301381, 11172158, 61303023, 11371184, 11426128), Science and Technology Development Fund of Tianjin Higher
Institutions (20121003), Doctoral Fund of Tianjin Normal University (52XB1202).

\end{document}